\documentclass[preprint,12pt]{elsarticle}

\journal{Springer Proceedings}
 
\usepackage{a4,amssymb,amsmath,amsthm}
\numberwithin{equation}{section}

\usepackage{hyperref}

\newtheorem*{theorem*}{Theorem}
\newtheorem{theorem}{Theorem}[section]
\newtheorem{remark}[theorem]{Remark}
\newtheorem*{proposition*}{Proposition}

\allowdisplaybreaks

\begin{document}

\begin{frontmatter}

\title{Symmetry-Preserving Finite-Difference Schemes and Auto-B\"acklund Transformations for the Schwarz Equation}

\author[label1]{Evgeniy I. Kaptsov\corref{cor1}}
\ead{evgkaptsov@math.sut.ac.th}
\cortext[cor1]{Corresponding author}
\affiliation[label1]{organization={School of Mathematics \& Geoinformatics, Institute of Science, Suranaree University of Technology},
city={Nakhon~Ratchasima},
postcode={30000}, 
country={Thailand}}

\author[label2]{Vladimir A. Dorodnitsyn}
\ead{Dorodnitsyn@Keldysh.ru}
\affiliation[label2]{organization={Keldysh Institute of Applied Mathematics, Russian Academy of Sciences},
addressline={Miusskaya~Pl.~4}, 
city={Moscow},
postcode={125047}, 
country={Russia}}

% use at least 70 and at most 150 words.
\begin{abstract}
It is demonstrated that one of the equations from the Lie classification list of second-order ODEs is a first integral of the Schwarz equation. 
As symmetry-preserving finite-difference schemes have been previously constructed for both equations, the preservation of a similar connection between these schemes is studied. It is shown that the schemes for the Schwarz equation and the second-order ODE can be related through a B\"acklund-type difference transformation. 

In addition, previously unexamined aspects of the difference scheme for the second-order ODE are discussed, including its singular solution and the complete set of difference first integrals.
\end{abstract}

\end{frontmatter}

\section{Inroduction}

Symmetries are closely related to the geometric properties of equations and play a crucial role in studying problems in mathematical physics~\cite{bk:Ovsiannikov1978,bk:Olver[1986]}. 
Knowledge of the symmetry of a differential equation often allows one to identify its conservation laws and exact solutions, as well as to reduce its order. 
As a particular case, for ordinary differential equations (ODEs), symmetries provide a means to obtain first integrals~(using Noether's theorem~\cite{art:Noether1918}, the adjoint equation method, or the direct method~\cite{bk:BlumanAnco2002}). In some cases, when a sufficient number of first integrals is known, an ODE can be completely integrated through algebraic calculations alone. Symmetry analysis also facilitates the reduction of higher-order equations to well-known reference equations of the same order or the lowering of their order by reduction on a subgroup. 
Among the numerous examples of such reductions, here we would like to highlight the works~\cite{bk:Tamizhmani_reductions2015,bk:Tamizhmani_Riccati2015} of Professor Kilkothur~M.~Tamizhmani and co-authors on the symmetries and reductions of third- and second-order ODEs to known cases, as well as the study of their Painlev\'e properties.

\medskip

In this paper, we consider two particular ODEs and their discrete symmetry-preserving analogues. 
The first ODE is the Schwarz equation
\begin{equation} \label{ODE3}
\frac{y^{\prime\prime\prime}}{y^\prime}
    - \frac{3}{2} \left(
        \frac{y^{\prime\prime}}{y^\prime}
    \right)^2 = 0,
\end{equation}
where the left-hand side represents the Schwarzian derivative. This derivative arises in various areas of mathematics, including classical complex analysis, one-dimensional dynamics, integrable systems, and conformal field theory~\cite{bk:OvsienkoTabachnikov_Schwarzian}. 

Equation (\ref{ODE3}) admits the following six-dimensional Lie algebra
\begin{equation} \label{alg6}
\def\arraystretch{2}
\begin{array}{ccc}
\displaystyle
X_1 = \frac{\partial}{\partial x},
\qquad
X_2 = x \frac{\partial}{\partial x},
\qquad
X_3 = x^2 \frac{\partial}{\partial x},
\\
\displaystyle
X_4 = \frac{\partial}{\partial y},
\qquad
X_5 = y \frac{\partial}{\partial y},
\qquad
X_6 = y^2 \frac{\partial}{\partial y}.
\end{array}
\end{equation}
In~\cite{bk:WintDorKapKozDAN2014,DorKapKozWin2015}, the adjoint equation method~\cite{bk:BlumanAnco2002} was employed to construct the first integrals of~(\ref{ODE3}), which will be discussed further. Through the integrals, one derives the general solution of~(\ref{ODE3}) as
\begin{equation} \label{solODE3}
y = \frac{1}{C_1 x + C_2} + C_3,
\end{equation}
where $C_1$, $C_2$, and $C_3$ are arbitrary constants.

\medskip

The second equation under consideration is
\begin{equation} \label{ODE2}
y^{\prime\prime} + 2\,\frac{y^\prime + C_0 {y^\prime}^{3/2} + {y^\prime}^2}{x - y} = 0,
\end{equation}
where $C_0$ is an arbitrary constant. 
Notice that this constant cannot be eliminated by any equivalence transformation, i.e., without altering the group structure of the equation.
This second-order ODE belongs to the Lie classification list~\cite{bk:Lie1888_classif,bk:HandbookLie_v1}. In his classification, Lie identified all second-order ODEs that admit nontrivial symmetries. 
Remarkably, most second-order ODEs whose solutions and integrals appear in standard handbooks are equivalent, up to point transformations, to equations from the Lie list. 

\smallskip

According to the classification, equation~(\ref{ODE2}) admits the tree-dimensional Lie algebra
\begin{equation} \label{alg3}
\displaystyle
Y_1 = \frac{\partial}{\partial x} + \frac{\partial}{\partial y},
\qquad
Y_2 = x \frac{\partial}{\partial x} + y \frac{\partial}{\partial y},
\qquad
Y_3 = x^2 \frac{\partial}{\partial x} + y^2 \frac{\partial}{\partial y}.
\end{equation}
Employing this symmetries and the Noether theorem, 
the general solution of~(\ref{ODE2}) is derived~\cite{bk:WintDorKapKozDAN2014} as
\begin{equation} \label{solODE2}
y = \frac{1}{A_0 (B_0 - A_0 x)} + \frac{B_0-C_0}{A_0},
\end{equation}
where $A_0$ and $B_0$ are constants of integration. 
Equation (\ref{ODE2}) also has the singular solution
\begin{equation} \label{specSolODE}
y = a_0 x + b_0, 
\end{equation}
where $a_0$ and $b_0$ are constants, and $a_0$ must satisfy the condition
\[
a_0(a_0  + C_0 \sqrt{a}_0 + 1) = 0.
\]

%TODO: Say something on discrete classification~\cite{bk:DorodnitsynKozlovWinternitz[2000]}
%
%Next, we focus on the symmetry-preserving discrete analogues of these equations and explore the relationships between them. 
%Finite-difference schemes that preserve all the symmetries of the original differential equations for both equations were previously constructed in~\cite{bk:Bourlioux2006,bk:DorKap2ndOrderODEScm2013}. 
%Such schemes are called invariant, and their construction theory has been well developed in recent decades~\cite{bk:Dorodnitsyn[2001],bk:LeviWintYamilov2023}. 
%In subsequent sections, we present invariant schemes for the ODEs under consideration and demonstrate that a connection can be established between these schemes using a B\"acklund-type finite-difference transformation.

\medskip

Notice that~(\ref{solODE2}) is a particular case of the solution~(\ref{solODE3}).
Indeed, a connection between~(\ref{ODE3}) and~(\ref{ODE2}) can be established.\footnote{For the case $C_0=0$, the connection between the solutions of these ODEs was noted long ago; see, e.g.,~\cite[p.~15]{bk:Ince1956ODEs}, which references the Mathematical Tripos, Part~I Examination, 1911.}
By solving~(\ref{ODE2}) for the constant~$C_0$, one obtains
\begin{equation} \label{ODE2C}
\displaystyle
\frac{(y - x) y^{\prime\prime} - 2 y^\prime \left(1 + y^\prime\right)}{2 {y^\prime}^{3/2}} = C_0.
\end{equation}
Differentiating (\ref{ODE2C}) with respect to $x$ yields
\begin{equation}  \label{ODE3_with_mult}
%D_x \left(
%    \frac{(y - x) y^{\prime\prime} - 2 y^\prime \left(1 + y^\prime\right)}{2 {y^\prime}^{3/2}}
%\right) =
\frac{y - x}{2 \sqrt{y^\prime}} \left(
    \frac{y^{\prime\prime\prime}}{y^\prime}
    - \frac{3}{2} \left(
        \frac{y^{\prime\prime}}{y^\prime}
    \right)^2
\right) 
= 0.
\end{equation}
This shows that~(\ref{ODE2C}) is a first integral of the Schwarz equation~(\ref{ODE3}), with the multiplier
\[
\frac{y - x}{2 \sqrt{y^\prime}}.
\]

\medskip

Similar to the Lie classification~\cite{bk:Lie1888_classif} for ODEs, a group classification of second-order finite-difference schemes was carried out in~\cite{bk:DorodnitsynKozlovWinternitz[2000]}. As a result, symmetry-preserving finite-difference schemes were constructed for the equations from the Lie list, including schemes for equation~(\ref{ODE2}). Such symmetry-preserving schemes are called \emph{invariant}, and their construction theory has been well developed in recent decades~\cite{bk:Dorodnitsyn[2011],bk:LeviWintYamilov2023}. 
Among invariant schemes, there exist \emph{exact} schemes, which are schemes whose solutions coincide exactly with the solutions of the corresponding differential equations at the nodes of the finite-difference mesh. Such schemes for ODEs were constructed, e.g., in~\cite{bk:Rodriguez_2004,bk:Kozlov_2007,bk:DorKap2ndOrderODEScm2013}. 
An exact scheme for equation~(\ref{ODE2}) was derived in~\cite{bk:DorKap2ndOrderODEScm2013}, while for the Schwarz equation~(\ref{ODE3}), it was obtained in~\cite{bk:Bourlioux2006,bk:Bourlioux2008}.

As discrete analogues for both equations~(\ref{ODE3}) and~(\ref{ODE2}) have been constructed, a natural question arises about the possibility of establishing a connection between them, similar to the one described above for the corresponding ODEs. 
In the following sections, we will demonstrate that there exists a finite-difference connection of a more complex nature, based on the application of B\"acklund-type finite-difference transformations.

\medskip

The remainder of the paper is organized as follows. 
In Section~{\ref{sec:Schemes}}, we provide a brief overview of the invariant finite-difference schemes for the Schwarz equation~(\ref{ODE3}) and the second-order ODE~(\ref{ODE2}). 
The connection between these schemes using B\"acklund-type finite-difference transformations is explored in Section~{\ref{sec:Backlund}}. 
Finally, Section~{\ref{sec:Concl}} discusses the implications of the established connection and presents conclusions based on the findings.

%%%%%%%%%%%%%%%%%%%%%%%%%%%%%%%%%%%%

\section{Invariant Finite-Difference Schemes for Eqs.~(\ref{ODE3}) and~(\ref{ODE2})}
\label{sec:Schemes}

We consider schemes on the following subset (a four-point stencil) of the finite-difference space
\begin{equation} \label{stencil_n}
\left(x_{n-1}, x_{n}, x_{n+1}, x_{n+2}, u_{n-1}, u_{n}, u_{n+1}, u_{n+2}\right)
\end{equation}
Further, the symbol $x_n=x(n)$ (or $t_n=t(n)$) is used as an approximation for the independent variable, and
$u_n=u(n)$ (or $y_n=y(n)$) is used as an approximation for the dependent variable at the point corresponding to the integer index~$n$ on the finite-difference grid. 
The index~$n$ is shifted left or right by the finite-difference shift operators~${\underset{-h}{S}}$ and~${\underset{+h}{S}}$, respectively.
The finite-difference derivatives are defined in terms of shifts as
\[
\underset{+h}{D} = \frac{{\underset{+h}{S}} - 1}{x_{n+1} - x_{n}}, 
\qquad
\underset{-h}{D} = \frac{1 - {\underset{-h}{S}}}{x_{n} - x_{n-1}}. 
\]
The operators ${\underset{+h}{S}}$ and ${\underset{-h}{S}}$ commute with the differentiations~$\underset{\pm h}{D}$ in any order. 

For brevity, we adopt Samarskii's notation~\cite{bk:SamarskyPopov_book[1992]} and write the variables~(\ref{stencil_n}) as
\begin{equation} \label{stencil}
\left(x_-, x, x_+, x_{++}, u_-, u, u_+, u_{++}\right).
\end{equation}
In general, for any function $f = f(n)$, we denote
\[
f_- = {\underset{-h}{S}}(f), \quad 
f_+ = {\underset{+h}{S}}(f), \quad
f_{++} = {\underset{+h}{S}}^2(f).
\]
We also denote the finite-difference derivatives and difference mesh steps as
\[
\displaystyle
u_{\bar{x}} = \underset{-h}{D}(u),
\qquad
u_x = \underset{+h}{D}(u),
\qquad
u_x^+ = \underset{+h}{D}(u_+) = \underset{+h}{S}(u_x),
\]
\[
h_- = x - x_-,
\qquad
h_+ = x_+ - x,
\qquad
h_{++} = x_{++} - x_+.
\]

In contrast to the continuous case~\cite{bk:Ovsiannikov1978,bk:Olver[1986]}, in variables (\ref{stencil}), the prolongation formula for the generator
\begin{equation} \label{generatorX}
\displaystyle
X = \xi \frac{\partial}{\partial x} + \eta \frac{\partial}{\partial u}
\end{equation}
is defined through finite-difference shifts~\cite{bk:Dorodnitsyn[2011]} as follows:
\[
\displaystyle
\underset{h}{\textrm{pr}} \, X = X 
+ \xi_- \frac{\partial}{\partial x_-} + \eta_- \frac{\partial}{\partial u_-}
+ \xi_+ \frac{\partial}{\partial x_+} + \eta_+ \frac{\partial}{\partial u_+}
+ \xi_{++} \frac{\partial}{\partial x_{++}} + \eta_{++} \frac{\partial}{\partial u_{++}}.
\]
A finite-difference equation $F(x_-, x, x_+, x_{++}, u_-, u, u_+, u_{++}) = 0$ is called \emph{invariant} with respect to the symmetry generator~(\ref{generatorX}) 
if 
\begin{equation} \label{FDinvCrit}
\underset{h}{\textrm{pr}} \, X (F) \big|_{[F]} = 0.
\end{equation}
Here $[F]$ indicates that (\ref{FDinvCrit}) is evaluated on solutions of~$F=0$ and all its finite-difference consequences.

An invariant finite-difference scheme on the stencil~(\ref{stencil}) is typically written as a system of two equations~\cite{bk:BakirovaDorodnitsynKozlov[2004],bk:Dorodnitsyn[2011]},
\begin{equation} \label{SchemeGen}
\def\arraystretch{1.5}
\begin{array}{l}
\Phi(x_-, x, x_+, x_{++}, u_-, u, u_+, u_{++}) = 0,
\\
\Omega(x_-, x, x_+, x_{++}, u_-, u, u_+, u_{++}) = 0,
\end{array}
\end{equation}
for which the invariance criteria holds:
\[
\underset{h}{\textrm{pr}} \, X_i (\Phi) \big|_{[\Phi],[\Omega]} = 0, 
\qquad
\underset{h}{\textrm{pr}} \, X_i (\Omega) \big|_{[\Phi],[\Omega]} = 0
\]
for all symmetries~$X_i$ admitted by the corresponding ODE. 
The first equation in~(\ref{SchemeGen}) serves as a finite-difference approximation of the ODE, while the second equation, referred to as the mesh equation, vanishes in the continuous limit, reducing to the trivial identity~$0 \equiv 0$. 
%Recall that an invariant scheme for an~ODE is called \emph{exact} if its solutions coincide exactly with the solutions of the~ODE at all nodes of the computational grid.

\medskip

An invariant scheme for~(\ref{ODE2}) was constructed by the authors in~\cite{bk:DorKap2ndOrderODEScm2013} 
and numerically investigated in~\cite{art:DorKapPreprint2014}. It can be written in the form
\begin{equation} \label{ScemeODE2}
\displaystyle
\theta \left( 
        \frac{\sqrt{u_{\bar{x}}}}{\sqrt{(x-u_-)(x_- - u)}}  
        - \frac{x_+ - u}{x_- - u}\frac{\sqrt{u_x}}{\sqrt{(x-u_+)(x_+ - u)}}  
    \right)
    + \frac{C \sqrt{1+\varepsilon} \, (x - x_-) u_{\bar{x}}}{(x- u_-)(x_- - u)} = 0,
\end{equation}
\begin{equation} \label{ScemeODE2Mesh}
\displaystyle
\frac{(x_+ - x) (u_+ - u)}{(x - u_+)(x_+ - u)} = \frac{(x - x_-) (u - u_-)}{(x - u_-)(x_- - u)} = \varepsilon,
\end{equation}
where $C$ is constant a corresponding to $C_0$, and the small parameter $\varepsilon$ characterizes the density of the finite-difference mesh, and $\theta$ ($|\theta| \to 1$) is a parameter. 
In case
\begin{equation} \label{theta}
\theta = \frac{\sqrt{1+\varepsilon}}{2}(|C| \sqrt{\varepsilon} - \sqrt{\varepsilon C^2 + 4}).
\end{equation}
the scheme becomes exact,\footnote{To simplify the calculations, in the next section we assume that~$\theta=1$. As verified through previous calculations carried out by the authors, this assumption does not qualitatively affect the subsequent results, while it substantially simplifies the equations.} 
i.e., its solutions at the grid nodes coincide with those of the corresponding ODE. 
As scheme (\ref{ScemeODE2}), (\ref{ScemeODE2Mesh}) is invariant, it admits the same algebra~(\ref{alg3}) as that of~(\ref{ODE2}).

By means of the difference analogue of the Noether theorem~\cite{bk:Dorodnitsyn[2011]}, 
the following three first integrals~\cite{bk:DorKap2ndOrderODEScm2013} were derived for the scheme 
\begin{equation} \label{schODE2_JJ}
\def\arraystretch{2}
\begin{array}{c}
\displaystyle
J_1 = \frac{C}{x_+ - u} - \frac{\theta (u_x + 1)}{\sqrt{u_x (x - u)(x_+ - u_+)}} = A,
\\
\displaystyle
J_2 = \frac{C x_+}{x_+ - u} - \frac{\theta (x_+ u_x + u)}{\sqrt{u_x (x - u)(x_+ - u_+)}} = B,
\\
\displaystyle
J_3 = \frac{(x_+ - x)(u_+ - u)}{(x - u_+)(x_+ - u)} = \varepsilon.
\end{array}
\end{equation}
Using these integrals, the scheme can be reduced to a~Riccati difference equation, which is solved by standard methods.

\medskip

For the Schwarz equation~(\ref{ODE3}), an invariant finite-difference scheme, known as the~Winternitz scheme, was proposed in~\cite{bk:Bourlioux2006,bk:Bourlioux2008} and has the form
\begin{equation} \label{WintScm}
\def\arraystretch{2.25}
\begin{array}{l}
\displaystyle
\frac{(y_{++} - y)(y_{+} - y_-)}{(y_{++} - y_+)(y - y_-)} - K = 0,
\\
\displaystyle
\frac{(t_{++} - t)(t_{+} - t_-)}{(t_{++} - t_+)(t - t_-)} - K = 0,
\end{array}
\end{equation}
where $K$ is constant. 
This scheme admits the same six-dimensional algebra~(\ref{alg6}) as the Schwarz equation.

In \cite{bk:WintDorKapKozDAN2014} and \cite{DorKapKozWin2015} scheme~(\ref{WintScm}) was integrated for arbitrary values of~$K$, as well as schemes of a somewhat more general form. This was first done for a difference equation of odd order using the difference analogue of the adjoint equation method~\cite{bk:BlumanAnco2002}, which is applicable even when the equation does not possess a Lagrangian or a Hamiltonian.

A comparison of the general solutions of scheme~(\ref{ScemeODE2}), (\ref{ScemeODE2Mesh}) and scheme~(\ref{WintScm}) 
found in~\cite{bk:DorKap2ndOrderODEScm2013,bk:WintDorKapKozDAN2014,DorKapKozWin2015} 
leads to the following relation between the constants~$K$ and~$C$
\[
K = \frac{\varepsilon C^2 + 4}{\varepsilon + 1}.
\]

In the case $K = 4$ ($C^2 = 4$), scheme~(\ref{WintScm}) is exact~\cite{bk:WintDorKapKozDAN2014}. 
For simplicity, further we restrict ourselves to this case, as the calculations for it are the most concise. The cases $K<4$ and $K>4$ can be addressed in a similar manner.

For $K=4$, the adjoint equation method gives six difference integrals of scheme~(\ref{WintScm}). 
Here, we present only the two integrals that we will need later, namely
\begin{equation}\label{WinScmInts}
\def\arraystretch{2.25}
\begin{array}{r}
\displaystyle
\frac{4}{y_+ - y_-} - \frac{1}{y_+ - y} - \frac{1}{y - y_-} = \textrm{const}, 
\\
\displaystyle
\frac{4}{t_+ - t_-} - \frac{1}{t_+ - t} - \frac{1}{t - t_-} = \textrm{const}.
\end{array}
\end{equation}
By means of the integrals, the general solution of scheme (\ref{WintScm}) is found:
\begin{equation}\label{K=4_gensol}
y_n = \frac{1}{c_1 n + c_2} + c_3,
\qquad
t_n = \frac{1}{c_4 n + c_5} + c_6,
\end{equation}
where $c_1$, $c_2$, ..., $c_6$ are constants of integration.

\smallskip

The exact solution of scheme~(\ref{ScemeODE2}), (\ref{ScemeODE2Mesh}) corresponding to the case $K=4$ is
\begin{equation} \label{C=2_gensol}
u_n = \frac{1}{A (B - A \, x_n)} + \frac{B-C}{A},
\qquad
x_n = \frac{\textrm{sgn} \, {C} \sqrt{1 + \varepsilon}}{A \sqrt{\varepsilon} (\rho + n)}
    + \frac{B - \textrm{sgn} \, {C}}{A},
\end{equation}
where $\rho$ is constant. Note that the constants $A$, $B$, and $C$ generally differ from the corresponding constants~$A_0$, $B_0$, and~$C_0$ in~(\ref{solODE2}). While it should be possible to adjust them to match, this is not our objective here.

\begin{remark}
Solution (\ref{C=2_gensol}) can be rewritten as
\[
\displaystyle 
u_n = \frac{\textsl{sgn} \, {C} \sqrt{\varepsilon}(\rho + n)}
    {A \left(\sqrt{\varepsilon} (\rho + n) - \sqrt{1+\varepsilon}\right)}
    + \frac{B - C}{A},
\quad
x_n = \frac{\textrm{sgn} \, {C} \sqrt{1 + \varepsilon}}{A \sqrt{\varepsilon} (\rho + n)}
    + \frac{B - \textrm{sgn} \, {C}}{A}.
\]
\end{remark}

\begin{remark}
Using~(\ref{C=2_gensol}), one can reconstruct the remaining first integral of~(\ref{ScemeODE2}), (\ref{ScemeODE2Mesh}) for $C^2=4$:
\[
%J_4 = \frac{1}{x_+ - u_+ + u - x} \left( \frac{\sqrt{(x_+ - u_+)(u_+ - u)(x - u)}}{\sqrt{x_+ - x}} + u - x \right) - n = \rho.
J_4 = \frac{1}{x_+  - x - u_+ + u} \left( u - x  + {\textsl{sgn} \, {C}} \sqrt{u_x (u - x)(u_+ - x_+)}\right) - (n + 1) = \rho.
\]
Apparently, this first integral cannot be obtained by means of the difference analogue of the Noether theorem due to its dependence on the index~$n$, 
which in some sense expresses the `non-locality' of the integral.  

%................
%
%
%We multiply by $x_+ - x = h_+$ to get
%\[
%\displaystyle
%h_+ J_4 = \frac{u - x  + {\textsl{sgn} \, {C}} \sqrt{u_x (u - x)(u_+ - x_+)}}{1 - u_x} - (n + 1) h_+.
%\]
%We notice that
%\[
%(n+1) h_+ = x_+ - x_0 = x - x_0 + h_+
%\]
%for some constant initial value $x_0$. Then, 
%\[
%h_+ J_4 = \frac{u - x  + {\textsl{sgn} \, {C}} \sqrt{u_x (u - x)(u_+ - x_+)}}{1 - u_x} - x + x_0 - h_+.
%\]
%Expanding this integral into a Taylor series gives:
%\[
%h_+ J_4 = \frac{u - x}{1 + \sqrt{u^\prime}} + x + O(h_+).
%\]
\end{remark}

\begin{remark}
As equation (\ref{ODE2}) has the singular solution~(\ref{specSolODE}), and scheme (\ref{ScemeODE2}), (\ref{ScemeODE2Mesh}) is exact, it should also possess a singular solution of a similar form. To verify this, substitute
\begin{equation} \label{specSolFD}
u_n = a x_n + b
\end{equation}
into (\ref{ScemeODE2Mesh}), resulting in the equation
\begin{multline*}
x_{n+1} = \frac{\varepsilon  }{a (2 \varepsilon + 1)} \left\{
        \left(
            (a^2 + 1) + \frac{a}{\varepsilon}
            + (a - 1) \sqrt{ a^2 + \frac{2 a (\varepsilon + 1)}{\varepsilon} + 1}
        \right) x_n
    \right.
    \\
    \left.
    + \, b \left(
        (a - 1)
        + \sqrt{ a^2 + \frac{2 a (\varepsilon + 1)}{\varepsilon} + 1}
    \right)
\right\}.
\end{multline*}
This is a linear equation solvable by standard methods. 
Substituting this and~(\ref{specSolFD}) into~(\ref{ScemeODE2}) taking into account (\ref{theta}) yields a quite cumbersome equation relating the values of~$\varepsilon$, $a$, and~$C$. This equation can be interpreted as a condition on the choice of the mesh density parameter~$\varepsilon$.
\end{remark}

\section{Establishing a Finite-Difference Connection Through B\"acklund-type Transformations (Case $K=C^2=4$)}
\label{sec:Backlund}

Suppose there exists a connection between scheme~(\ref{ScemeODE2}), (\ref{ScemeODE2Mesh}) and scheme (\ref{WintScm}) through finite-difference differentiation, similar to the corresponding ODEs. To verify this, equation~(\ref{ScemeODE2}) is first rewritten in an integral form by solving for~$C$:
\begin{equation} \label{ScmODE2Int1}
\frac{(x- u_-)(x_- - u)}{(u - u_-) \sqrt{1+\varepsilon}}\left( 
        \frac{x_+ - u}{x_- - u}\frac{\sqrt{u_x}}{\sqrt{(x-u_+)(x_+ - u)}}   
        - \frac{\sqrt{u_{\bar{x}}}}{\sqrt{(x-u_-)(x_- - u)}} 
    \right)
    = C.
\end{equation}
Next, assuming that (\ref{ScmODE2Int1}) can be reduced to a first integral of the Winternitz scheme~(\ref{WintScm}), and given that~(\ref{WintScm}) is symmetric with respect to~$y$ and~$t$, there should be another integral obtained from~(\ref{ScmODE2Int1}) by interchanging~$u$ and~$x$. 
As the calculations confirm, the following equation holds on the mesh~(\ref{ScemeODE2Mesh})
\begin{equation} \label{ScmODE2Int2}
\frac{(u- x_-)(u_- - x)}{(x - x_-)\sqrt{1+\varepsilon}}\left( 
        \frac{u_+ - x}{u_- - x}\frac{\sqrt{x_u}}{\sqrt{(u-x_+)(u_+ - x)}}  
        - \frac{\sqrt{x_{\bar{u}}}}{\sqrt{(u-x_-)(u_- - x)}}  
    \right)
    = \widetilde{C},
\end{equation}
where we have denoted
\[
x_u = \frac{x_+ - x}{u_+ - u},
\qquad
x_{\bar{u}} = \frac{x - x_-}{u - u_-},
\qquad
\widetilde{C} = \frac{\sqrt{\varepsilon} - \sqrt{1+\varepsilon}}{\sqrt{\varepsilon} + \sqrt{1+\varepsilon}} \, C.
\]
%Note that 
%\begin{equation}
%J_5 = C / \widetilde{C} = \frac{(x_+ - u)(x - x_-)}{(x_+ - x)(u - x_-)} 
%= \frac{\sqrt{1+\varepsilon} - \sqrt{\varepsilon}}{\sqrt{1+\varepsilon} + \sqrt{\varepsilon}} = \textrm{const}
%\end{equation}
%is a first integral of (\ref{ScmODE2Int1}), (\ref{ScmODE2Int2}).

Thus, scheme (\ref{ScmODE2Int1}), (\ref{ScmODE2Int2}) has been obtained in integral form, and it is equivalent to scheme~(\ref{ScemeODE2}), (\ref{ScemeODE2Mesh}). 

Differentiating the scheme by applying the operator~$\underset{+h}{D}$ to (\ref{ScmODE2Int1}) and (\ref{ScmODE2Int2}), one derives the four-point scheme
\begin{equation}\label{ScmODE3fromScmODE2}
\begin{array}{c}
\displaystyle
\frac{(x_+ - u)(x - u_+)}{u_+ - u}\left( 
        \frac{\sqrt{u_x}}{\sqrt{(x_+-u)(x - u_+)}}  
        - \frac{x_{++} - u_+}{x - u_{+}}\frac{\sqrt{u_x^+}}{\sqrt{(x_+-u_{++})(x_{++} - u_+)}}  
    \right) {}
\\ \displaystyle
{} - \frac{(x- u_-)(x_- - u)}{u - u_-}\left( 
        \frac{\sqrt{u_{\bar{x}}}}{\sqrt{(x-u_-)(x_- - u)}}  
        - \frac{x_+ - u}{x_- - u}\frac{\sqrt{u_x}}{\sqrt{(x-u_+)(x_+ - u)}}  
    \right) = 0,
\\ {}
\\ \displaystyle
\frac{(u_+ - x)(u - x_+)}{x_+ - x}\left( 
        \frac{\sqrt{x_u}}{\sqrt{(u_+-x)(u - x_+)}}  
        - \frac{u_{++} - x_+}{u - x_{+}}\frac{\sqrt{x_u^+}}{\sqrt{(u_+-x_{++})(u_{++} - x_+)}}  
    \right) {}
\\ \displaystyle
{} - \frac{(u- x_-)(u_- - x)}{x - x_-}\left( 
        \frac{\sqrt{x_{\bar{u}}}}{\sqrt{(u-x_-)(u_- - x)}}  
        - \frac{u_+ - x}{u_- - x}\frac{\sqrt{x_u}}{\sqrt{(u-x_+)(u_+ - x)}}  
    \right) = 0.
\end{array}
\end{equation}
The first equation of (\ref{ScmODE3fromScmODE2}) approximates the Schwarz equation, while the second one serves as a mesh equation. 
Clearly, (\ref{ScmODE2Int1}) and (\ref{ScmODE2Int2}) are first integrals of (\ref{ScmODE3fromScmODE2}), as system~(\ref{ScmODE3fromScmODE2}) 
was obtained by differentiating them.

It might be expected that scheme~(\ref{ScmODE3fromScmODE2}) could be transformed into the Winternitz scheme (\ref{WintScm}) via some point transformation; however, this is not the case. It can be verified that scheme~(\ref{ScmODE3fromScmODE2}) admits only the three-dimensional subalgebra~(\ref{alg3}) of the six-dimensional algebra~(\ref{alg6}), which is admitted by the Winternitz scheme. This implies that scheme~(\ref{ScmODE3fromScmODE2}) cannot be transformed into scheme~(\ref{WintScm}) by a point transformation, and therefore, more general transformations must be considered.  
Note also that scheme~(\ref{ScmODE3fromScmODE2}) remains exact on the subset~(\ref{C=2_gensol}) of the solutions of the Winternitz scheme.

\smallskip

To proceed, consider B\"acklund-type transformations, starting with the differential case. 
Recall that the equation 
\begin{equation} \label{Backl_gen}
\mathcal{B}\left(x, u, y, u^\prime, y^\prime, u^{\prime\prime}, y^{\prime\prime}\right) = 0
\end{equation}
defines a B\"acklund transformation that relates equations
\begin{equation} \label{ODE3B1}
\frac{u^{\prime\prime\prime}}{u^\prime}
    - \frac{3}{2} \left(
        \frac{u^{\prime\prime}}{u^\prime}
    \right)^2 = 0
\end{equation}
and
\begin{equation} \label{ODE3B2}
\frac{y^{\prime\prime\prime}}{y^\prime}
    - \frac{3}{2} \left(
        \frac{y^{\prime\prime}}{y^\prime}
    \right)^2 = 0
\end{equation}
if the compatibility of (\ref{Backl_gen}) and (\ref{ODE3B1}) implies (\ref{ODE3B2}), and the compatibility of (\ref{Backl_gen}) and~(\ref{ODE3B2}) implies~(\ref{ODE3B1}). 
In this particular case, both~(\ref{ODE3B1}) and~(\ref{ODE3B2}) are Schwarz equations, and a transformation of this type is referred to as an auto-B\"acklund transformation. 
By analyzing the compatibility conditions of (\ref{Backl_gen}), (\ref{ODE3B1}), and~(\ref{ODE3B2}), 
it can be demonstrated that an equation of the form
\begin{equation}   \label{Backl_gen1}
\displaystyle
\mathcal{B}\left( 
    \frac{u^{\prime\prime}}{{u^{\prime}}^{3/2}}, 
    u - \frac{2 {u^\prime}^2}{u^{\prime\prime}}, 
    \frac{y^{\prime\prime}}{{y^{\prime}}^{3/2}}, 
    y - \frac{2 {y^\prime}^2}{y^{\prime\prime}}
\right) = 0,
\end{equation}
where the arguments of $\mathcal{B}$ are first integrals of~(\ref{ODE3}), defines an auto-B\"acklund transformation.  
As a straightforward example of such a transformation, one can consider
\begin{equation} 
\mathcal{B} = u - \frac{2 {u^\prime}^2}{u^{\prime\prime}} + \alpha \left(\frac{y^{\prime\prime}}{{y^{\prime}}^{3/2}}\right) = 0,
\end{equation}
where $\alpha$ is some nonzero constant which can be expressed in terms of integration constants of~(\ref{ODE3B1}) and~(\ref{ODE3B2}).

\smallskip

Assuming that scheme~(\ref{ScmODE3fromScmODE2}) corresponds to equation~(\ref{ODE3B1}) and scheme~(\ref{WintScm}) corresponds to equation~(\ref{ODE3B2}), and using the known first integrals of these schemes, one can construct a difference analogue of the transformation~(\ref{Backl_gen1}). It is important to note that in the difference case, two equations, $\mathcal{B}_1$ and $\mathcal{B}_2$, are needed to connect independent first integrals from both schemes, as the systems~(\ref{ScmODE3fromScmODE2}) and~(\ref{WintScm}) include mesh equations. An example of such a transformation is
\begin{equation} \label{Backl}
\def\arraystretch{2.5}
\begin{array}{c}
\displaystyle
\mathcal{B}_1 = \frac{(x- u_-)(x_- - u)}{u - u_-}\left( 
        \frac{\sqrt{u_{\bar{x}}}}{\sqrt{(x-u_-)(x_- - u)}}  
        - \frac{x_+ - u}{x_- - u}\frac{\sqrt{u_x}}{\sqrt{(x-u_+)(x_+ - u)}}  
    \right) 
    \\ \displaystyle {} + \alpha_1 \left( \frac{4}{y_+ - y_-} - \frac{1}{y_+ - y} - \frac{1}{y - y_-} \right) = 0,
\\ {}
\displaystyle
\mathcal{B}_2 = \frac{(u- x_-)(u_- - x)}{x - x_-}\left( 
        \frac{\sqrt{x_{\bar{u}}}}{\sqrt{(u-x_-)(u_- - x)}}  
        - \frac{u_+ - x}{u_- - x}\frac{\sqrt{x_u}}{\sqrt{(u-x_+)(u_+ - x)}}  
    \right) 
    \\ \displaystyle {} + \alpha_2 \left( \frac{4}{t_+ - t_-} - \frac{1}{t_+ - t} - \frac{1}{t - t_-} \right) = 0,
\end{array}
\end{equation}
where $\alpha_1$ and $\alpha_2$ are nonzero constants which can be expressed in terms of the integration constants of the schemes.
Here, four integrals (\ref{WinScmInts}), (\ref{ScmODE2Int1}), and (\ref{ScmODE2Int2}) are employed, though this transformation is not the only possible one. 
It is easy to verify that system~(\ref{ScmODE3fromScmODE2}) follows from the compatibility conditions of~(\ref{WintScm}) and~(\ref{Backl}), 
and system~(\ref{WintScm}) follows from the compatibility conditions of~(\ref{ScmODE3fromScmODE2}) and~(\ref{Backl}). In other words, 
a B\"acklund-type difference transformation relating schemes (\ref{WintScm}) and~(\ref{ScmODE3fromScmODE2}) was constructed.

\medskip

Based on the above construction, one formulates the following more general
\begin{proposition*}
Given two schemes for ODEs of the same order $k>1$ with at least two known first integrals each, a finite-difference B\"acklund-type transformation can be constructed utilizing these integrals to relate the schemes. 
If both schemes approximate the same ODE, this transformation reduces to an auto-B\"acklund transformation in the continuous limit. 
\end{proposition*}

It is worth noting that discrete B\"acklund-type transformations have previously been applied to differential-difference equations~\cite{bk:LeviWintYamilov2023} and to problems related to the linearization of finite-difference schemes~\cite{bk:DorodnLinearization2006}.

\begin{remark}
As the general solutions~(\ref{K=4_gensol}) and~(\ref{C=2_gensol}) of schemes~(\ref{ScmODE3fromScmODE2}) and~(\ref{WintScm}) are known, substituting these solutions,
\begin{equation}
x_{\pm} = x_{\pm}(n), \quad
t_{\pm} = t_{\pm}(n), \quad
u_{\pm} = u_{\pm}(n), \quad
y_{\pm} = y_{\pm}(n), \quad \cdots,
\end{equation}
allows explicit relations between $x$, $u$, $t$, and $y$ to be expressed in terms of the index $n$ and integration constants. 
The expressions for transformation~(\ref{Backl}) are quite cumbersome and are not presented here.
\end{remark}

\section{Conclusions}
\label{sec:Concl}

An example was considered in which the second-order ODE~(\ref{ODE2}) is transformed into a form where it becomes a first integral~(\ref{ODE2C}) of the Schwarz equation~(\ref{ODE3}). It was shown that, in the discrete case, this structure does not hold for the known finite-difference schemes.

By differentiating the exact scheme for the second-order equation, written in terms of the first integrals~(\ref{ScmODE2Int1}) and~(\ref{ScmODE2Int2}), a scheme for equation~(\ref{ODE3}) is obtained that does not coincide with the known exact scheme~(\ref{WintScm}). Instead, the scheme~(\ref{ScmODE3fromScmODE2}) is derived, which is exact only for a subset of solutions of~(\ref{ScmODE2Int1}), (\ref{ScmODE2Int2}) and admits only the subalgebra~(\ref{alg3}) of the six-dimensional algebra~(\ref{alg6}). 
This implies that the two schemes for the Schwarz equation cannot be related by a point transformation.

Nevertheless, a connection between schemes~(\ref{ScmODE3fromScmODE2}) and~(\ref{WintScm}) can be established through a finite-difference transformation~(\ref{Backl}) of the Bäcklund type, utilizing the known first integrals of these two schemes.
In the continuous limit, both schemes~(\ref{WintScm}) and~(\ref{ScmODE3fromScmODE2}) reduce to the Schwarz equation (up to some multipliers), while the discrete transformation~(\ref{Backl}) simplifies to an auto-B\"acklund transformation, mapping the Schwarz equation onto itself.

The results suggest that multiple exact schemes for the same ODE may exist, connected by point transformations~(for isomorphic Lie algebras) 
or more complex transformations like discrete B\"acklund-type transformations.

\section*{Acknowledgements}

This research was partially supported by Suranaree University of Technology and the NSRF via the Program Management 
Unit for Human Resources \& Institutional Development, Research and Innovation (PMU-B) (Grant~B13F660067). 
%E.I.K. also expresses gratitude to Woxsen University, Hyderabad, for its warm welcome during the initial preparation of the manuscript.
E.I.K. is also grateful to Prof.~S.~V.~Meleshko for support and fruitful discussions during the current postdoctoral fellowship~(Full-time66/06/2023).

%\cite{bk:Bourlioux2006}
%\cite{bk:Bourlioux2008}
%\cite{bk:WintDorKapKozDAN2014}
%\cite{DorKapKozWin2015}
%\cite{bk:LeviWintYamilov2023}
%\cite{bk:DorodnLinearization2006}

%
% ---- Bibliography ----
%

%\bibliographystyle{unsrt}
%\bibliography{references}

\end{document}